\documentclass[12pt]{article}
\usepackage{xcolor} %%%%%package for colors
\usepackage[utf8]{inputenc}
\usepackage{amsmath}
\usepackage{amssymb}
\usepackage{mathrsfs}
\usepackage{amsthm}
\usepackage{cite}
\usepackage{hyperref}
\usepackage{enumerate}
\usepackage[left=0.7in,right=0.7in,top=0.7in,bottom=0.7in]{geometry}
\usepackage{titlesec}
\titleformat*{\section}{\large\bfseries}
\newtheorem{theorem}{Theorem}[section]
\newtheorem{lemma}[theorem]{Lemma}
\newtheorem{corollary}[theorem]{Corollary}

\newtheorem{remark}[theorem]{Remark}
\numberwithin{equation}{section}

\title{Self-adjoint and co-isometry composition and weighted composition operator on general weighted Hardy space}
\author{\large Anuradha Gupta and Geeta Yadav$^*$}
\date{}
\begin{document}
\maketitle
\begin{abstract} In this paper we study the self-adjoint and co-isometry composition and weighted composition operator on $H_E(\zeta).$ We also discuss the conditions under which adjoint operator of a composition and weighted composition operator on $H_E(\zeta)$ to be some composition and weighted composition operator, respectively.


\end{abstract}    
\section{Introduction and Preliminaries}
Consider a sequence of positive real numbers $\zeta=(\zeta_n)_{n\in \mathbb{N}_0}$ such that
\begin{equation} \label{condition for entire functions}
lim_{n\rightarrow \infty} (\zeta_n)^{1/n}= \infty
\end{equation}
where $\mathbb{N}_0$ denotes the set $\mathbb{N} \cup \{0\} $. The weighted Hardy space of entire functions on complex plane induced by $\zeta,$ denoted by $H_E(\zeta)$ (see \cite{Main Paper on H_E(beta)}), is defined as
$$H_E(\zeta)=\Big\{g(z)=\sum_{n=0}^{\infty} b_n z^n : ||g||^2=\sum_{n=0}^{\infty} |b_n|^2 {\zeta_n}^2 <\infty \Big\}.$$ 
It is a Hilbert space with the inner product defined by 
\begin{equation} \label{inner product define}
\langle g,h \rangle=\sum_{n=0}^{\infty} b_n \overline{c_n} {\zeta_n}^2 
\end{equation}
for $g(z)=\sum_{n=0}^{\infty} b_n z^n$ and $h(z)=\sum_{n=0}^{\infty} c_n z^n$ in $H_E(\zeta).$ Functions in Hilbert space $H_E(\zeta)$ are entire if and only if \eqref{condition for entire functions} is satisfied \cite{Hilbert space of entire fun}.

For space $H_E(\zeta),$  the sequence $(e_n)_{n\in \mathbb{N}_0}$ where $e_n(z)=\frac{1}{\zeta_n}z^n$ for $n\in \mathbb{N}_0,$ forms an orthonormal basis. 
The well known Fock space, $\mathcal{F}^2(\mathbb{C}),$ is a particular case of $H_E(\zeta)$ when sequence $\zeta=(\sqrt{n!})_{n\in \mathbb{N}_0}.$\\

Let $X$ be a space whose elements are analytic functions defined on some domain $V$ in $\mathbb{C}.$ For analytic functions  $\Upsilon:V \longrightarrow \mathbb{C}$  and $\Phi: V \longrightarrow V,$ the weighted composition operator is defined by
$$C_{\Upsilon,\Phi}f=\Upsilon \cdot( f \circ \Phi)$$
for $f\in X.$ If $\Upsilon \equiv 1,$ then it becomes the composition operator $C_{\Phi}f=( f \circ \Phi)$ induced by $\Phi.$ \\

Le \cite{Normal and isometric} completely characterized the bounded weighted composition operator on $\mathcal{F}^2(\mathbb{C}).$ Many researchers have studied the self-adjoint and co-isometric weighted composition operators on Hilbert spaces. 
Zhao and Pang \cite{Zhao and Pang} and Zhao \cite{Normal weig comp op on the fock sp} characterized the self-adjoint weighted composition operators on $\mathcal{F}^2(\mathbb{C})$ and $\mathcal{F}^2(\mathbb{C}^N),$ respectively. In 2016, Zhao \cite{A note on inv weigh comp on the fock sp}  obtained that weighted composition operator is co-isometry on $\mathcal{F}^2(\mathbb{C}^N)$ if and only if it is unitary.

In this paper we study the relation between the function-theoretic behaviours of $\Phi$ and $ \Upsilon$ with the operator-theoretic properties of $C_{\Upsilon,\Phi}.$ This paper is organized as follows:\\
In section 2, we discuss the necessary and sufficient conditions for the adjoint operator of a composition and weighted composition operator on $H_E(\zeta)$ to be some composition and weighted composition operator, respectively. Further, with the help of these results we discuss self-adjoint composition and weighted composition operator on $H_E(\zeta).$ The last section deals with co-isometry composition and weighted composition operator on $H_E(\zeta).$
  
%%%%%%%%%%%%%%%%%%%%%%%%%%%%%%%%%%%%%%%%%%%%%%%%%%%%%%%%%%%%%%%%%%%%%%%%%%%%%%%%%%%%%%%%
\section{Self-adjoint composition and weighted composition operators}
The Hilbert space $H_E(\zeta)$ is a reproducing kernel Hilbert space (RKHS). For each $p\in \mathbb{C},$ there exists a function $K_{p},$ known as reproducing kernel, defined by
\begin{equation} \label{reproducing kernel function}
K_p(w)= \sum_{n=0}^{\infty} \frac{1}{{\zeta_n}^2} \bar{p}^n w^n  \,\,\,\, \text{for all}\,\, w\in \mathbb{C}
\end{equation} 
such that $f(p)=\langle f, K_p \rangle $ for each function $f$ in $H_E(\zeta).$ In particular, for $\mathcal{F}^2(\mathbb{C}), K_p(w)=e^{\bar{p} w}.$
First, we will discuss some well known properties of reproducing kernel which will be used in the subsequent results. 
From equation \eqref{reproducing kernel function}, it follows that $K_0(w)=\frac{1}{{\zeta_0}^2}$ and $K_p(0)=\frac{1}{{\zeta_0}^2}.$ Also,
\begin{equation} \label{property of conjugate reproducing kernel}
\overline{K_p(w)}= \overline{\langle K_{p}, K_{w} \rangle}=\langle K_{w}, K_{p} \rangle=K_w(p) \,\,\,\,\text{for all} \,\,  p,w \in \mathbb{C}
\end{equation}
and
\begin{equation} \label{norm of reproducing kernel}
||K_{p}||^2=\langle K_{p}, K_{p} \rangle= K_{p}(p) \,\,\,\,\text{for all} \,\,  p \in \mathbb{C}.
\end{equation}
For a bounded weighted composition operator $C_{\Upsilon,\Phi}$ on $H_E(\zeta)$
\begin{equation*}
\langle g,C_{\Upsilon,\Phi}^* K_{z} \rangle=\langle C_{\Upsilon,\Phi} g, K_{z} \rangle= \Upsilon(z) \langle g,K_{\Phi(z)} \rangle= \langle g,\overline{\Upsilon(z)} K_{\Phi(z)} \rangle \,\,\text{for all} \,\, g\in H_E(\zeta),\, z \in \mathbb{C}.
\end{equation*}
Thus, $C_{\Upsilon,\Phi}^* K_{z}=\overline{\Upsilon(z)} K_{\Phi(z)}$ for all $z\in \mathbb{C}.$\\

%%%%%%%%%%%%%%%%%%%%%%%%%%%%%%%%%%%%%%%%%%%%%%%%%%%%%%%%%%%%%%%%%%%%%%%%%%%%%%%%%%%%%%
In the following result $\Phi_{1}$ and $\Phi_{2}$  (for $N=1$) are taken according to the necessary condition obtained by 
Doan et al. (\cite{Com Op on several variable}, Theorem 4.6) for $C_{\Phi_1}$ and $C_{\Phi_2}$ to be bounded on $H_E(\zeta),$ respectively. Tan and Khoi \cite{Main Paper on H_E(beta)} studied the boundedness of composition operator on $H_E(\zeta)$ and obtained many useful results.
%%%%%%%%%%%%%%%%%%%%%%%%%%%%%%%%%%%%%%%%%%%%%%%%%%%%%%%%%%%%%%%%%%%%%%%%%%
The necessary and sufficient condition for the adjoint of a bounded composition operator to be some composition operator on $H_E(\zeta)$ is as follows: 

\begin{theorem} \label{thm necessary condition on adjoint com op}
Let $\Phi_{1}(z)=\mu z+c$ and $\Phi_{2}(z)=\nu z+d$ with $|\mu| \leq 1$ and $|\nu| \leq 1,$ respectively,  such that $C_{\Phi_{1}}$ and $C_{\Phi_{2}}$ are bounded on $H_E(\zeta).$ Then, the adjoint of a composition operator $C_{\Phi_{1}}$ is composition operator $C_{\Phi_{2}}$ if and only if
$$c=0,\, d=0 \,\, \text{and} \,\, \bar{\mu}=\nu.$$
\end{theorem}
\textbf{Proof.} Let $C_{\Phi_{1}}^*=C_{\Phi_{2}}.$ Then, for all $z,w\in \mathbb{C}$
\begin{align}
(C_{\Phi_{1}}^*K_{z})(w)&=(C_{\Phi_{2}}K_{z})(w) \notag\\
 K_{\Phi_{1}(z)}(w)&=K_{z}(\Phi_{2}(w)). \label{adjoint of comp is wei adjoint 1}
\end{align}
Taking $z=0,$ for all $w\in \mathbb{C}$ we obtain   
\begin{align*}
K_{\Phi_{1}(0)}(w) &= \frac{1}{{\zeta_0}^2} \\
\sum_{n=0}^{\infty} \frac{1}{{\zeta_n}^2} \, {\bar{c}}^n w^n &=\frac{1}{{\zeta_0}^2}.
\end{align*}
Comparing both sides we get $c=0$ because ${\zeta_n}>0,$ for all $n \in \mathbb{N}_0.$ Now taking $w=0$ in equation \eqref{adjoint of comp is wei adjoint 1} we get
\begin{align*} 
 K_{\Phi_{1}(z)}(0)&=K_{z}(\Phi_{2}(0)) \\
\frac{1}{{\zeta_0}^2} &= \sum_{n=0}^{\infty} \frac{1}{{\zeta_n}^2} \, {\bar{z}}^n {d}^n.
\end{align*}
Comparing both sides we get $d=0.$ Since $c=0$ and $d=0$ in $\Phi_{1}$ and $\Phi_{2},$ respectively, therefore, from equation \eqref{adjoint of comp is wei adjoint 1} we get
$$\sum_{n=0}^{\infty} \frac{1}{{\zeta_n}^2} \, {\bar{\mu}}^n {\bar{z}}^n w^n= \sum_{n=0}^{\infty} \frac{1}{{\zeta_n}^2} \, {\bar{z}}^n {\nu}^n  w^n \,\,\, \text{for all} \,\, z,w\in \mathbb{C}.$$
Comparing both sides we get $\bar{\mu}=\nu$ because the above equation holds for all (nonzero) $z, w \in \mathbb{C}.$ \\
Conversely, let $\Phi_{1}(z)=\mu z,$ $\Phi_{2}(z)=\nu z$ and $\bar{\mu}=\nu.$ Since $C_{\Phi_{1}}$ and $C_{\Phi_{2}}$ are bounded on $H_E(\zeta)$ and set $Span\{K_{z} : z\in \mathbb{C}\}$ is dense in $H_E(\zeta),$ so to prove $C_{\Phi_{1}}^*=C_{\Phi_{2}}$ on $H_E(\zeta),$  it is sufficient to prove that $C_{\Phi_{1}}^*K_{z}=C_{\Phi_{2}}K_{z}$ for all $z\in \mathbb{C}.$ Now, for all $z,w \in \mathbb{C}$
\begin{align*}
(C_{\Phi_{1}}^*K_{z})(w) &=K_{\Phi_{1}(z)}(w) \\
                         &=\sum_{n=0}^{\infty} \frac{1}{{\zeta_n}^2} \, {\bar{\mu}}^n {\bar{z}}^n w^n \\
                         &=\sum_{n=0}^{\infty} \frac{1}{{\zeta_n}^2} \, {\bar{z}}^n {\nu}^n  w^n \\
                         &=(C_{\Phi_{2}}K_{z})(w).
\end{align*}
Hence, $C_{\Phi_{1}}^*=C_{\Phi_{2}}.$ \\
%%%%%%%%%%%%%%%%%%%%%%%%%%%%%%%%%%%%%%%%%%%%%%%%%%%%%%%%%%%%%%%%%%%%%%%%%%%%%

In above Theorem if $\Phi_{1}$ and $\Phi_{2}$ are same, then we obtain the following characterization for self-adjoint composition operator on $H_E(\zeta).$

\begin{corollary} \label{self adoint mu z function}
Let $\Phi(z)=\nu z+d$ for all $z\in \mathbb{C}$ with $|\nu| \leq 1$ such that $C_{\Phi}$ is bounded on $H_E(\zeta).$  Then, $C_{\Phi}$ is self-adjoint if and only if $d=0 \,\, \text{and} \,\, \nu \,\, \text{is a real number}.$
\end{corollary}
%%%%%%%%%%%%%%%%%%%%%%%%%%%%%%%%%%%%%%%%%%
Zhao and Pang \cite{Zhao and Pang} obtained the necessary and sufficient condition for the adjoint operator of a bounded weighted composition operator on $\mathcal{F}^2(\mathbb{C})$  to be some weighted composition operator.

\begin{theorem} \label{thm necessary condition on adjoint weig com op}
Let $\Phi_{1}(z)=\mu z+c$ and $\Phi_{2}(z)=\nu z+d$ for $c, d\in \mathbb{C}$ with $|\mu| \leq 1$ and $|\nu| \leq 1,$ respectively, and $\Upsilon_{1},\Upsilon_{2}$ be entire functions on $\mathbb{C}$ such that $C_{\Upsilon_{1},\Phi_{1}}$ and $C_{\Upsilon_{2},\Phi_{2}}$ are bounded on $H_E(\zeta).$ Then, the adjoint of a weighted composition operator $C_{\Upsilon_{1},\Phi_{1}}$ is weighted composition operator $C_{\Upsilon_{2},\Phi_{2}}$ if and only if
\begin{equation} \label{first condition for adjoint of weighted is weighted}
\Upsilon_{1}(z) ={\zeta_0}^2 \Upsilon_{1}(0) \, K_{\Phi_{2}(0)}(z) \,\, \text{and} \,\,\, \Upsilon_{2}(z)={\zeta_0}^2 \, \overline{\Upsilon_{1}(0)} K_{\Phi_{1}(0)}(z) \,\, \text{for all}\,\, z\in \mathbb{C}
\end{equation}
and one of the following is true
\begin{itemize}
\item[(i)] If $\Upsilon_{1}(0)=0$ then  $\Upsilon_{1}\equiv 0$ and $\Upsilon_{2}\equiv 0$ consequently $C_{\Upsilon_{1},\Phi_{1}} \equiv 0$ and $C_{\Upsilon_{2},\Phi_{2}} \equiv 0.$ \\
\item[(ii)] If $\Upsilon_{1}(0) \neq 0$ then 
\begin{equation} \label{adjoint of weghted com op}
K_{z}(\Phi_{2}(0)) K_{\Phi_{1}(z)}(w)= K_{\Phi_{1}(0)}(w) K_{z}(\Phi_{2}(w))  \,\,\,\, \text{for all}\,\, z,w\in \mathbb{C}. 
\end{equation}
\end{itemize}
\end{theorem}
\textbf{Proof.} Let $C_{\Upsilon_{1},\Phi_{1}}^*=C_{\Upsilon_{2},\Phi_{2}}.$ Then, for all $z,w\in \mathbb{C}$
\begin{align}
(C_{\Upsilon_{1},\Phi_{1}}^*K_{z})(w)&=(C_{\Upsilon_{2},\Phi_{2}}K_{z})(w) \notag\\
\overline{\Upsilon_{1}(z)} K_{\Phi_{1}(z)}(w)&=\Upsilon_{2}(w)K_{z}(\Phi_{2}(w)). \label{adjoint of weigh comp is wei adjoint 1}
\end{align}
Taking $w=0$ we obtain  
\begin{align*}
\overline{\Upsilon_{1}(z)} &={\zeta_0}^2  \Upsilon_{2}(0)K_{z}(\Phi_{2}(0)) \,\,\,\, \text{for all}\, z\in \mathbb{C}. 
\end{align*}
Taking $z=0$ we get that $\overline{\Upsilon_{1}(0)} =\Upsilon_{2}(0).$ Thus, 
\begin{equation}\label{adjoint of weigh comp is wei adjoint 2}
\overline{\Upsilon_{1}(z)} ={\zeta_0}^2 \overline{\Upsilon_{1}(0)}\, K_{z}(\Phi_{2}(0)) \,\,\,\,  \text{for all}\,\, z\in \mathbb{C}.
\end{equation}
Taking conjugate on both sides we get $\Upsilon_{1}(z) ={\zeta_0}^2 \Upsilon_{1}(0) \, K_{\Phi_{2}(0)}(z).$
Combining equations (\ref{adjoint of weigh comp is wei adjoint 1}) and (\ref{adjoint of weigh comp is wei adjoint 2}) we obtain 
\begin{equation} \label{adjoint of weigh comp is wei adjoint 3}
{\zeta_0}^2  \overline{\Upsilon_{1}(0)}\, K_{z}(\Phi_{2}(0)) K_{\Phi_{1}(z)}(w)=\Upsilon_{2}(w)K_{z}(\Phi_{2}(w)) \,\,\,\,  \text{for all}\,\, z,w\in \mathbb{C}.
\end{equation}
Taking $z=0$ we get %in equation (\ref{adjoint of weigh comp is wei adjoint 3}) 
\begin{align} \label{adjoint of weigh comp is wei adjoint 4}
{\zeta_0}^2 \overline{\Upsilon_{1}(0)} \frac{1}{{\zeta_0}^2}    K_{\Phi_{1}(0)}(w)&=\Upsilon_{2}(w) \frac{1}{{\zeta_0}^2}\,\,\,\,  \text{for all}\,\, w\in \mathbb{C} \notag \\
{\zeta_0}^2 \overline{\Upsilon_{1}(0)}    K_{\Phi_{1}(0)}(w)&=\Upsilon_{2}(w) \,\,\,\,  \text{for all}\,\, w\in \mathbb{C}.
\end{align}
Combining equations (\ref{adjoint of weigh comp is wei adjoint 3}) and (\ref{adjoint of weigh comp is wei adjoint 4}) 
\begin{equation} \label{adjoint of weigh comp is wei adjoint 5}
\overline{\Upsilon_{1}(0)} K_{z}(\Phi_{2}(0)) K_{\Phi_{1}(z)}(w)=\overline{\Upsilon_{1}(0)} K_{\Phi_{1}(0)}(w) K_{z}(\Phi_{2}(w))  \,\,\,\, \text{for all}\,\, z,w\in \mathbb{C}.
\end{equation}
Case I. If $\Upsilon_{1}(0)=0$ then from equations (\ref{adjoint of weigh comp is wei adjoint 2}) and (\ref{adjoint of weigh comp is wei adjoint 4}) it follows that $\Upsilon_{1}\equiv 0$ and $\Upsilon_{2}\equiv 0.$\\
Case II. If $\Upsilon_{1}(0)\neq 0$ then from equation (\ref{adjoint of weigh comp is wei adjoint 5}) we obtain
\begin{equation} \label{adjoint of weigh comp is wei adjoint 6}
 K_{z}(\Phi_{2}(0)) K_{\Phi_{1}(z)}(w)= K_{\Phi_{1}(0)}(w) K_{z}(\Phi_{2}(w))  \,\,\,\, \text{for all}\,\, z,w\in \mathbb{C}.
\end{equation}
Conversely, let $\Upsilon_{1}$ and $\Upsilon_{2}$ be as given in equation \eqref{first condition for adjoint of weighted is weighted}. If $\Upsilon_{1}(0)=0$ then trivially $C_{\Upsilon_{1},\Phi_{1}}^*=C_{\Upsilon_{2},\Phi_{2}}$ so assume that $\Upsilon_{1}(0)\neq 0$ and equation \eqref{adjoint of weghted com op} holds. For all $z,w \in \mathbb{C}$
\begin{align*}
(C_{\Upsilon_{1},\Phi_{1}}^*K_{z})(w)&=\overline{\Upsilon_{1}(z)} K_{\Phi_{1}(z)}(w)  \\
             &={\zeta_0}^2 \overline{\Upsilon_{1}(0)}\, K_{z}(\Phi_{2}(0)) K_{\Phi_{1}(z)}(w). 
 \end{align*}
From equation \eqref{adjoint of weghted com op} it follows that, for all $z,w \in \mathbb{C}$        
 \begin{align*}            
(C_{\Upsilon_{1},\Phi_{1}}^*K_{z})(w) &={\zeta_0}^2 \overline{\Upsilon_{1}(0)} \, K_{\Phi_{1}(0)}(w) K_{z}(\Phi_{2}(w)) \\
             &=\Upsilon_{2}(w) K_{z}(\Phi_{2}(w)) \\
             &=(C_{\Upsilon_{2},\Phi_{2}}K_{z})(w) .
\end{align*}
Thus, $C_{\Upsilon_{1},\Phi_{1}}^*K_{z}=C_{\Upsilon_{2},\Phi_{2}}K_{z}$ for all $z\in \mathbb{C}$ and hence, $C_{\Upsilon_{1},\Phi_{1}}^*=C_{\Upsilon_{2},\Phi_{2}}.$\\

%%%%%%%%%%%%%%%%%%%%%%%%%%%%%%%%%%%%%%%%%%%%%%%%%
In above theorem, we are unable to get more precise condition for $\Phi_{1}$ and $\Phi_{2}$ in Case II on $H_E(\zeta).$

%%%%%%%%%%%%%%%%%%%%%%%%%%%%%%%%%%%%%%%%%%%%%%%%%%%%%%%%%%%
\begin{theorem} \label{thm general case self adjoint  nu z}
Let $\Phi(z)=\nu z$  with $|\nu| \leq 1$ and $\Upsilon$ be entire functions on $\mathbb{C}$ such that $C_{\Upsilon,\Phi}$ is bounded on $H_E(\zeta).$ Then, $C_{\Upsilon,\Phi}$ is self-adjoint if and only if
$$ \Upsilon(z)=\overline{\Upsilon(0)} $$
and one of the following is true
\begin{itemize}
\item[(i)] For $\Upsilon(0)=0,$ $\Upsilon\equiv 0$ consequently $C_{\Upsilon,\Phi} \equiv 0.$ \\
\item[(ii)] For $\Upsilon(0)\neq 0,$  $\nu$ is a real number.
\end{itemize}
\end{theorem}
\textbf{Proof.} Clearly, $\Phi(0)=0.$ Let $C_{\Upsilon,\Phi}$ be self-adjoint on $H_E(\zeta).$ Then, by Theorem \ref{thm necessary condition on adjoint weig com op} for $\Phi_{1}=\Phi=\Phi_{2}$ and $\Upsilon_{1}=\Upsilon=\Upsilon_{2},$ equation \eqref{first condition for adjoint of weighted is weighted}  gives 
\begin{equation*} 
\Upsilon(z)=\overline{\Upsilon(0)} \,\,\,\, \text{for all} \,\, z\in \mathbb{C}.
\end{equation*} 
Case I. If $\Upsilon(0)=0,$ then $\Upsilon\equiv 0.$\\
Case II. If $\Upsilon(0) \neq 0,$  then from equation \eqref{adjoint of weghted com op} we obtain
\begin{equation}  
\frac{1}{{\zeta_0}^2}K_{\Phi(z)}(w)= \frac{1}{{\zeta_0}^2}K_{z}(\Phi(w)) \,\,\,\, \text{for all} \, z, w\in \mathbb{C}.
\end{equation}
Taking $w=z,$ for all $z\in \mathbb{C}$ we get
\begin{align}
K_{\Phi(z)}(z)&= K_{z}(\Phi(z)) \\
\sum_{n=0}^{\infty} \frac{1}{{\zeta_n}^2} (\bar{\nu} \bar{z})^n z^n &=\sum_{n=0}^{\infty} \frac{1}{{\zeta_n}^2} \overline{z}^n (\nu z)^n \\
\sum_{n=0}^{\infty} \frac{1}{{\zeta_n}^2} {\bar{\nu}}^n |z|^{2n} &=\sum_{n=0}^{\infty} \frac{1}{{\zeta_n}^2}  {\nu} ^n |z|^{2n}.
\end{align}
This implies $\bar{\nu}=\nu$ because the above equation holds for all $z (\neq 0) \in \mathbb{C}.$\\
Conversely, let $\Upsilon(z)=\overline{\Upsilon(0)}$ then taking $z=0$ we get $\Upsilon(0)=\overline{\Upsilon(0)}.$ If $\Upsilon(0)=0,$ then $C_{\Upsilon,\Phi} \equiv 0$ is trivially self-adjoint. So assume that $\Upsilon(0)\neq 0$ and $\bar{\nu}=\nu.$ Now
\begin{equation*}
(C_{\Upsilon,\Phi}^*K_{z})(w)=\overline{\Upsilon(z)} K_{\Phi(z)}(w)=\Upsilon(0) \sum_{n=0}^{\infty} \frac{1}{{\zeta_n}^2} (\bar{\nu} \bar{z})^n w^n \,\,\,\text{for all}\,\, z,w\in \mathbb{C}.
\end{equation*}
Since $\bar{\nu}=\nu$ and $\Upsilon(0)=\overline{\Upsilon(0)},$ therefore, for all $z,w\in \mathbb{C}$  
\begin{equation*}
(C_{\Upsilon,\Phi}^*K_{z})(w)=\overline{\Upsilon(0)} \sum_{n=0}^{\infty} \frac{1}{{\zeta_n}^2} (\overline{z})^n (\nu w)^n=\Upsilon(w) K_{z}(\Phi(w))=(C_{\Upsilon,\Phi}K_{z})(w).
\end{equation*}
Thus, $C_{\Upsilon,\Phi}^*K_{z}=C_{\Upsilon,\Phi}K_{z}$ for all $z\in \mathbb{C}$ and hence, $C_{\Upsilon,\Phi}^*=C_{\Upsilon,\Phi}.$

\begin{theorem} \label{thm general case self adjoint  nu z+d}
Let $\Phi(z)=\nu z+d$  for $\nu \neq 0, d\neq 0$ with $|\nu| \leq 1$ and $\Upsilon$ be entire functions on $\mathbb{C}$ such that $C_{\Upsilon,\Phi}$ is bounded on $H_E(\zeta).$ If $C_{\Upsilon,\Phi}$ is self-adjoint, then
$$ \Upsilon(z)=\alpha  K_{\Phi(0)}(z)\,\, \text{where} \,\, \alpha={\zeta_0}^2 \overline{\Upsilon(0)} \,\, \text{is a real number} $$
and one of the following is true
\begin{itemize}
\item[(i)] For $\Upsilon(0)=0,$ $\Upsilon\equiv 0$ consequently $C_{\Upsilon,\Phi} \equiv 0.$ \\
\item[(ii)] For $\Upsilon(0)\neq 0,$  $\nu$ is a real number.
\end{itemize}
\end{theorem}
\textbf{Proof.} Let $C_{\Upsilon,\Phi}$ be self-adjoint on $H_E(\zeta).$  Then, by Theorem \ref{thm necessary condition on adjoint weig com op} for $\Phi_{1}=\Phi=\Phi_{2}$ and $\Upsilon_{1}=\Upsilon=\Upsilon_{2},$ equation \eqref{first condition for adjoint of weighted is weighted}  gives  
\begin{align*}
 \Upsilon(z)&={\zeta_0}^2 \overline{\Upsilon(0)} K_{\Phi(0)}(z).  
\end{align*}
Taking $z=0,$ we get $\Upsilon(0)= \overline{\Upsilon(0)}$. Since $\zeta_0>0,$ therefore,
\begin{equation} \label{self adjoint nu z+d equ 2*}
 \Upsilon(z)=\alpha  K_{\Phi(0)}(z)\,\,\text{for all}\,\, z\in \mathbb{C}\,\, \text{where} \,\, \alpha={\zeta_0}^2 \overline{\Upsilon(0)} \,\, \text{and} \,\, \bar{\alpha}=\alpha.
\end{equation}
Case I. If $\Upsilon(0)=0$ then, $\Upsilon\equiv 0.$\\
Case II. If $\Upsilon(0)\neq 0,$ then from equation \eqref{adjoint of weghted com op} we obtain 
\begin{align}
  K_{z}(\Phi(0)) K_{\Phi(z)}(w)&=  K_{\Phi(0)}(w) K_{z}(\Phi(w)) \,\,\,\, \text{for all}\,\, z,w\in \mathbb{C}. \notag
\end{align} 
Taking $z=\big(\frac{-d}{\nu} \big),$ $w=\big(\frac{-d}{\nu} \big)$ and using $\Phi(0)=d,$ $\Phi\big(\frac{-d}{\nu} \big)=0$ we get
\begin{align}
 %%K_{\big(\frac{-d}{\nu}\big)}(d) K_{\Phi\left(\frac{-d}{\nu}\right)}(w)&=  K_{d} \big(\frac{-d}{\nu}\big) K_{z}(\Phi(\big(\frac{-d}{\nu})) 
\sum_{n=0}^{\infty} \frac{1}{{\zeta_n}^2} {\Big(\overline{\frac{-d}{\nu}}\Big)}^n {d\,}^n &=\sum_{n=0}^{\infty} \frac{1}{{\zeta_n}^2} \bar{d\,}^n {\Big(\frac{-d}{\nu} \Big)}^n \notag \\
\sum_{n=0}^{\infty} \frac{1}{{\zeta_n}^2} (-1)^n \, \frac{|d|^{2n}} {\bar{\nu}^n} &=\sum_{n=0}^{\infty} \frac{1}{{\zeta_n}^2} (-1)^n \, \frac{|d|^{2n}} {{\nu}^n}.  \notag
\end{align} 
This implies $\bar{\nu}=\nu$ because $d \neq 0.$
%%%%%%%%%%%%%%%%%%%%%%%%%%%%%%%%%%%%%%%%%%%%%%%%%%%%%%%%%%%%%%%%%%%%%%%%%%%%%%%%%%%
\begin{corollary} \label{self adjoint constant case}
Let $\Phi(z)=d$ for all $z\in \mathbb{C}$ and $\Upsilon$ be entire functions on $\mathbb{C}$ such that $C_{\Upsilon,\Phi}$ is bounded on $H_E(\zeta).$ Then, $C_{\Upsilon,\Phi}$ is self-adjoint if and only if
 $$ \Upsilon(z)=\alpha  K_{\Phi(0)}(z) \,\,\, \text{for all}\,\, z \in\, \mathbb{C}$$ 
 where $\alpha={\zeta_0}^2 \overline{\Upsilon(0)}$ is a real number.
%$$\Upsilon(w)={\zeta_0}^2\overline{\Upsilon(0)} K_{\Phi(0)}(w) \,\,\, \forall w \in\, \mathbb{C}.$$  
\end{corollary}
\textbf{Proof.} 
 Let  $C_{\Upsilon,\Phi}$ be self-adjoint. Then, $\Upsilon$ follows exactly on the same lines as in Theorem \ref{thm general case self adjoint  nu z+d}.\\  
Conversely, let $\Upsilon(z)=\alpha  K_{\Phi(0)}(z).$ Since $\alpha={\zeta_0}^2 \overline{\Upsilon(0)}$ is real and $\zeta_0 >0,$ therefore, $\Upsilon(0)=\overline{\Upsilon(0)}.$ Now consider, for all $z,w\in \mathbb{C}$
\begin{align*}
(C_{\Upsilon,\Phi}^*K_{z})(w)&=\overline{\Upsilon(z)} K_{\Phi(z)}(w) \\
                         &={\zeta_0}^2 \Upsilon(0) K_{z}(\Phi(0)) K_{\Phi(z)}(w). 
\end{align*}
From the fact that $\overline{\Upsilon(0)}=\Upsilon(0),$ $\Phi(z)=d$ and in particular, $\Phi(0)=d,$ we get 
\begin{align*}
(C_{\Upsilon,\Phi}^*K_{z})(w)&={\zeta_0}^2 \overline{\Upsilon(0)} K_{z}(d) K_{d}(w) \\
                         &=  {\zeta_0}^2 \overline{\Upsilon(0)} K_{\Phi(0)}(w) K_{z}(\Phi(w))\\ 
                         &=\Upsilon(w)K_{z}(\Phi(w))\\
                         &=(C_{\Upsilon,\Phi}K_{z})(w).
\end{align*}
Thus, $C_{\Upsilon,\Phi}^*=C_{\Upsilon,\Phi}$ and consequently, $C_{\Upsilon,\Phi}$ is self adjoint.\\
%%%%%%%%%%%%%%%%%%%%%%%%%%%%%%%%%%%%%%%%%%%%%%%%%%%%%%%%%%%%%%%%%%%%%%%%%%%%%%%%%%%%%%%%%%%

In above theorem $\nu=0,$ a real number in $\Phi(z)=\nu z+d$. The following result follows from  Theorem \ref{thm general case self adjoint  nu z}, Theorem \ref{thm general case self adjoint  nu z+d} and Corollary \ref{self adjoint constant case}.
\begin{theorem} \label{ general case involving all case}
Let $\Phi(z)=\nu z+d$ with $|\nu| \leq 1$ and $\Upsilon$ be entire functions on $\mathbb{C}$ such that $C_{\Upsilon,\Phi}$ is bounded on $H_E(\zeta).$ If $C_{\Upsilon,\Phi}$ is self-adjoint then
$$ \Upsilon(z)=\alpha  K_{\Phi(0)}(z)\,\, \text{where} \,\, \alpha={\zeta_0}^2 \overline{\Upsilon(0)} \,\, \text{is a real number} $$
and one of the following is true
\begin{itemize}
\item[(i)] For $\Upsilon(0)=0,$ $\Upsilon\equiv 0$ consequently $C_{\Upsilon,\Phi} \equiv 0.$ 
\item[(ii)] For $\Upsilon(0)\neq 0,$  $\nu$ is a real number.
\end{itemize}
\end{theorem}

%%%%%%%%%%%%%%%%%%%%%%%%%%%%%%%%%%%%%%%%%%%%%%%%%%%%%%%%%%%%%%%%%%%%%%%%%%%%%%%%%%%%%%%%%%%
%%%%%%%%%%%%%%%%%%%%%%%%%%%%%%%%%%%%%%%%%%%%%%%%%%%%%%%
\section{Co-isometry and adjoint of composition and weighted composition operators}
In 2014, Le \cite{Normal and isometric} completely characterized the isometric weighted composition operator on $\mathcal{F}^2(\mathbb{C}).$ It is well known that a bounded linear operator $S$ on a Hilbert space $H$ is co-isometry if $S^*$ is isometry, equivalently, $SS^*=I_{H}.$ The following result characterizes the co-isometry composition operator.

\begin{theorem}  \label{composition op co-isometry}
Let $\Phi(z)=\nu z-d$ for all $z\in \mathbb{C}$ with $|\nu| \leq 1$ such that $C_{\Phi}$ is bounded on $H_E(\zeta).$  Then, $C_{\Phi}$ is co-isometry if and only if $d=0 \,\, \text{and} \,\, |\nu|= 1.$
\end{theorem}
\textbf{Proof.} Let $C_{\Phi}$ be co-isometry on $H_E(\zeta).$ Now, for all $z,w\in \mathbb{C}$ we have
\begin{align}
(C_{\Phi} C_{\Phi}^*K_{z})(w)&=K_{z}(w) \notag \\
(C_{\Phi} K_{\Phi(z)})(w)&=K_{z}(w) \notag \\
K_{\Phi(z)}(\Phi(w))&=K_{z}(w). \label{Com Coisometry eq 1} 
\end{align}
This implies that $\nu \neq 0$ because for $\nu=0,$  $\Phi \equiv -d$ and then, $C_{\Phi}$ cannot be co-isometry.
Taking $w=\big(\frac{d}{\nu}\big)$ in equation \eqref{Com Coisometry eq 1} and using the fact that $\Phi\big(\frac{d}{\nu}\big)=0$ we get
\begin{align*}
 \frac{1}{{\zeta_0}^2}&=\sum_{n=0}^{\infty} \frac{1}{{\zeta_n}^2} \bar{z\,}^n \Big(\frac{d}{\nu}\Big)^n \,\,\,\text{for all} \,\,z\in \mathbb{C} \\
0&=\sum_{n=1}^{\infty} \frac{1}{{\zeta_n}^2} \Big(\frac{d}{\nu}\Big)^n \bar{z\,}^n. 
\end{align*}
Since ${\zeta_n} >0$ for all $n\in \mathbb{N}$ and above equation holds for all $z (\neq 0) \in \mathbb{C},$ comparing coefficients on both sides we get $d=0$ this implies $\Phi(z)=\nu z.$  Further, taking $w=z$ in equation \eqref{Com Coisometry eq 1}, for all $z\in \mathbb{C}$ we obtain
\begin{align*}
K_{\Phi(z)}(\Phi(z))&=K_{z}(z)  \\
\sum_{n=0}^{\infty} \frac{1}{{\zeta_n}^2} (\bar{\nu} \bar{z})^n (\nu z)^n &= \sum_{n=0}^{\infty} \frac{1}{{\zeta_n}^2} {\bar{z\,}}^n  {z}^n  \\
\sum_{n=0}^{\infty} \frac{1}{{\zeta_n}^2} |\nu|^{2n} |z|^{2n} &= \sum_{n=0}^{\infty} \frac{1}{{\zeta_n}^2} |z|^{2n}. 
\end{align*}
Compairing coefficients on both sides we get $|\nu|=1.$\\
Conversely, let $\Phi(z)=\nu z-d$ with $d=0$ and $|\nu|=1$ this implies $\Phi(z)=\nu z.$ Now, for all $z,w \in \mathbb{C}$
\begin{align}
(C_{\Phi} C_{\Phi}^*K_{z})(w)&=K_{\Phi(z)}(\Phi(w)) \notag \\
   &=\sum_{n=0}^{\infty} \frac{1}{{\zeta_n}^2} |\nu|^{2n} {\bar{z}}^n  w^n \\
   &=\sum_{n=0}^{\infty} \frac{1}{{\zeta_n}^2} {\bar{z}}^n  w^n.
\end{align}
Thus $C_{\Phi} C_{\Phi}^*K_{z}=K_{z}$ for all $z\in \mathbb{C}$ and hence, $C_{\Phi}$ is co-isometry.\\
%%%%%%%%%%%%%%%%%%%%%%%%%%%%%%%%%%%%%%%%%%%%%%%%%%%%%%%%%%%%%%%%%%%%%%%%%%%%%%%%%%%%%%%%%%%%%%

It can be easily checked that for $\Phi(z)=\nu z$ where $|\nu|=1,$
$\langle f \circ \Phi ,g \circ \Phi \rangle=\langle f  ,g \rangle \,\,\, \text{for all} \,\, f,g\in H_E(\zeta).$
This implies $\langle C_{\Phi}f,C_{\Phi}g \rangle= \langle f  ,g \rangle$ and hence, $C_{\Phi}$ is isometry. This along with Theorem \ref{composition op co-isometry} leads to the following result.

%%%%%%%%%%%%%%%%%%%%%%%%%%%%%%%%%%%%%%%%%%%%%%%%%%%%%%%%%%%%%%%%%%%%%%%%%%%%%%%%
\begin{corollary} 
Let $\Phi(z)=\nu z-d$ for all $z\in \mathbb{C}$ with $|\nu| \leq 1$ such that $C_{\Phi}$ is bounded on $H_E(\zeta).$  Then, $C_{\Phi}$ is co-isometry if and only if it is unitary.
\end{corollary}

%%%%%%%%%%%%%%%%%%%%%%%%%%%%%%%%%%%%%%%%%%%%%%%%%%%%%%%%%%%%%%%%%%%%%%%%%%%%%
\begin{lemma} \label{thm Entire fu space coisometry  Phi equal nu z}
Let $\Phi(z)=\nu z$ for all $z \in \mathbb{C}$ with $|\nu| \leq 1$ and $\Upsilon$ be entire function on $\mathbb{C}$ such that $C_{\Upsilon,\Phi}$ is bounded on $H_E(\zeta).$ Then, $C_{\Upsilon,\Phi}$ is co-isometry if and only if $$|\nu|=1\,\, \text{and} \,\, |\Upsilon(z)|=1 \,\, \text{where}\,\,\Upsilon(z)=\frac{1}{\overline{\Upsilon(0)}} \,\,\text{for all}\,\, z\in \mathbb{C}.$$ 
\end{lemma}
\textbf{Proof.} Let $C_{\Upsilon,\Phi}$ be co-isometry on $H_E(\zeta).$ Then, for all $z,w\in \mathbb{C}$ we have
\begin{align}
(C_{\Upsilon,\Phi} C_{\Upsilon,\Phi}^*K_{z})(w)&=K_{z}(w) \notag \\
\Upsilon(w)(( \overline{\Upsilon(z)} K_{\Phi(z)})\circ \Phi)(w)&=K_{z}(w) \notag \\
\Upsilon(w) \overline{\Upsilon(z)} K_{\Phi(z)}(\Phi(w))&=K_{z}(w). \label{coisometry nu z eq 1}
\end{align}
Taking $z=0$ we get
\begin{align} \label{coisometry nu z eq 1*}
\Upsilon(w) \overline{\Upsilon(0)} K_{\Phi(0)}(\Phi(w))&=\frac{1} {{\zeta_0}^2} \,\,\,\, \text{for all}\,\, w\in \mathbb{C}.    \notag 
\end{align}
This implies, $\Upsilon(0) \neq 0$ because if $\Upsilon(0)= 0$ then left side of equation will be equal to zero but $\zeta_0> 0$ so the right side is nonzero. Applying the fact that $\Phi(0)=0$ we get 
\begin{align}
\Upsilon(w)  &=\frac{1} {\overline{\Upsilon(0)}} \,\,\,\, \text{for all}\,\, w\in \mathbb{C}. 
\end{align}
Taking $w=0,$ we get $\Upsilon(0)=\dfrac{1} {\overline{\Upsilon(0)}},$ i.e., $|\Upsilon(0)|=1.$ From equation (\ref{coisometry nu z eq 1*}) we get
\begin{equation}   \label{coisometry nu z eq 2}
|\Upsilon(w)| =1 \,\,\, \text{for all}\,\, w\in\mathbb{C}. 
\end{equation}
Taking $w=z$ in equation (\ref{coisometry nu z eq 1}) and using equation  (\ref{coisometry nu z eq 2}), for all $z\in \mathbb{C}$ we get
\begin{align*}
\Upsilon(z) \overline{\Upsilon(z)} K_{\Phi(z)}(\Phi(z))&=K_{z}(z) \\
|\Upsilon(z)|^2 K_{\Phi(z)}(\Phi(z))&=K_{z}(z) \\
\sum_{n=0}^{\infty} \frac{1}{{\zeta_n}^2} (\bar{\nu} \bar{z})^n (\nu z)^n &=\sum_{n=0}^{\infty} \frac{1}{{\zeta_n}^2} \bar{z\,}^n z^n  \\
\sum_{n=0}^{\infty} \frac{1}{{\zeta_n}^2} {|\nu|}^{2n} {|z|}^{2n}  &=\sum_{n=0}^{\infty} \frac{1}{{\zeta_n}^2} {|z|}^{2n}.
\end{align*}
Comparing coefficients on both sides we get $|\nu|=1.$\\
Conversely, let $|\nu|=1\,\, \text{and} \,\, |\Upsilon(z)|=1 \,\, \text{for}\,\,\Upsilon(z)=\dfrac{1}{\overline{\Upsilon(0)}}.$ Now, for all $z,w \in \mathbb{C}$
\begin{equation}
(C_{\Upsilon,\Phi} C_{\Upsilon,\Phi}^*K_{z})(w)=\Upsilon(w) \overline{\Upsilon(z)} K_{\Phi(z)}(\Phi(w))=\frac{1}{|\Upsilon(0)|^2} \sum_{n=0}^{\infty} \frac{1}{{\zeta_n}^2} (\bar{\nu} \bar{z})^n (\nu w)^n. 
\end{equation}
From the fact that $|\Upsilon(0)|=1$ and $|\nu|=1,$ we get 
\begin{equation}
(C_{\Upsilon,\Phi} C_{\Upsilon,\Phi}^*K_{z})(w)= \sum_{n=0}^{\infty} \frac{1}{{\zeta_n}^2} {|\nu|}^{2n} \, {\bar{z\,} }^n {w }^{n}=\sum_{n=0}^{\infty} \frac{1}{{\zeta_n}^2} \, {\bar{z\,} }^n {w }^{n}=K_{z}(w) \,\,\, \text{for all}\,\, z, w\in \mathbb{C}.
\end{equation}
Thus, $C_{\Upsilon,\Phi} C_{\Upsilon,\Phi}^*K_{z}=K_{z}$ for all $z\in \mathbb{C}$ and hence $C_{\Upsilon,\Phi} C_{\Upsilon,\Phi}^*=I_{H_E(\zeta)}.$\\
%%%%%%%%%%%%%%%%%%%%%%%%%%%%%%%%%%%%%%%%%%%%%%%%%%%%%%%%%%%%%%%%%%%%%%%%%%%%%%%%%%%%%%%%%%%%%%%%%

From basic properties of adjoint operators it easily follows that if $S$ is an unitary operator on a Hilbert space and scalar $|\alpha|=1$ then $\alpha S$ is unitary. Since $C_{\Phi}$ is unitary  when $\Phi(z)=\nu z$ for $|\nu|= 1$ and $C_{\Upsilon,\Phi}=\alpha C_{\Phi}$ (in Lemma \ref{thm Entire fu space coisometry  Phi equal nu z}) where $\alpha=\dfrac{1}{\overline{\Upsilon(0)}}$ with $|\alpha|=1,$  therefore, $C_{\Upsilon,\Phi}$ is unitary and we have the following result: 

%%%%%%%%%%%%%%%%%%%%%%%%%%%%%%%%%%%%%%%%%%%%%%%%%%%%%%%%%%%%%%%%%%%%%%%%%%%%%%%%%%%%%%%%%%%%%
\begin{corollary} 
Let $\Phi(z)=\nu z$ for all $z \in \mathbb{C}$ with $|\nu| \leq 1$ and $\Upsilon$ be entire function on $\mathbb{C}$ such that $C_{\Upsilon,\Phi}$ is bounded on $H_E(\zeta).$ Then, $C_{\Upsilon,\Phi}$ is co-isometry if and only if it is unitary.
\end{corollary}

%%%%%%%%%%%%%%%%%%%%%%%%%%%%%%%%%%%%%%%%%%%%%%%%%%%%%%%%%%%%%%%%%%%%%%%%%%%%%%%%%%%%%%%%%
\begin{remark}
We can easily check that if $\Phi$ is a constant function say $\Phi \equiv c$ then for any entire function  $\Upsilon$  on $\mathbb{C},$ bounded operator $C_{\Upsilon,\Phi}$ cannot be co-isometry on $H_E(\zeta).$ Infact, the necessary condition for such $C_{\Upsilon,\Phi}$ to be co-isometry is that 
$$\Upsilon(z)=\frac{1}{{\zeta_0}^2 \overline{\Upsilon(0)} K_{c}(c)} \,\,\,\,\,\, \text{for all}\,\,\, z\in \mathbb{C}$$
that is, $\Upsilon$ is a constant function but for constant functions $\Phi$ and $\Upsilon,$ $C_{\Upsilon,\Phi}C_{\Upsilon,\Phi}^* \neq I_{H_E(\zeta)}$ and consequently $C_{\Upsilon,\Phi}$ cannot be co-isometry. 
\end{remark}

In view of this discussion, in the following result it is clear that for co-isometry $\nu$ is nonzero.

%%%%%%%%%%%%%%%%%%%%%%%%%%%%%%%%%%%%%%%%%%%%%%%%%%
\begin{theorem} \label{co-isometry mu z-d function}
Let $\Phi(z)=\nu z-d$ with $|\nu| \leq 1$ and $\Upsilon$ be entire functions on $\mathbb{C}$ such that $C_{\Upsilon,\Phi}$ is bounded on $H_E(\zeta).$  If $C_{\Upsilon,\Phi}$ is co-isometry then 
$$ \Upsilon(z)=\alpha \frac {K_{\bar{\nu} d} (z)} { ||K_{\bar{\nu} d}||} \, \,\,\, \text{for all} \,\, z\in \mathbb{C} $$
where $\alpha=\dfrac{1}{\overline{\Upsilon(0)} ||K_{\bar{\nu} d}||},$  $|\alpha|=\zeta_0$ and $|\nu|=1.$
\end{theorem}
\textbf{Proof.} Let $C_{\Upsilon,\Phi}$ be co-isometry on $H_E(\zeta).$ Then, for all $z,w\in \mathbb{C}$ we have
\begin{align}
(C_{\Upsilon,\Phi} C_{\Upsilon,\Phi}^*K_{z})(w)&=K_{z}(w) \notag \\
\Upsilon(w) \overline{\Upsilon(z)} K_{\Phi(z)}(\Phi(w))&=K_{z}(w). \label{coisometry nu z-d eq 1}
\end{align} 
Taking $w=0$ we get
\begin{align} \label{coisometry nu z-d eq 2}
\Upsilon(0) \overline{\Upsilon(z)} K_{\Phi(z)}(\Phi(0))&=\frac{1}{{\zeta_0}^2} \,\,\, \text{for all}\, \, z \in \mathbb{C}.
\end{align}
This implies that $\Upsilon(z) \neq 0$ because if $\Upsilon(z)= 0,$ then left side will be zero but right side is nonzero. From equation (\ref{coisometry nu z-d eq 2}) using the fact that $\Phi(0)=-d$ we get 
\begin{equation} \label{coisometry nu z-d eq 3}
\overline{\Upsilon(z)}= \frac{1}{{\zeta_0}^2 \Upsilon(0) K_{\Phi(z)}(-d)} \,\,\, \text{for all}\, \, z \in \mathbb{C}.
\end{equation}
Taking $z=0$ we get
\begin{equation} \label{coisometry nu z-d eq 4}
\overline{\Upsilon(0)}= \frac{1}{{\zeta_0}^2 \Upsilon(0) K_{(-d)}(-d)} \,\,\,\text{or}\,\, \,\, |\Upsilon(0)|^2=\frac{1}{{\zeta_0}^2  K_{(-d)}(-d)}  \,\, \,\,\,\text{for all}\, \, z \in \mathbb{C} .
\end{equation}
Taking $w=(\frac{d}{\nu})$ in equation (\ref{coisometry nu z-d eq 1}) and using the fact that $\Phi(\frac{d}{\nu})=0$ we get 
\begin{align} 
\Upsilon\Big(\frac{d}{\nu}\Big)  \overline{\Upsilon(z)} \frac{1}{{\zeta_0}^2} &=K_{z}\Big(\frac{d}{\nu}\Big) \notag\\
\overline{\Upsilon(z)} &={\zeta_0}^2 \frac{K_{z}\Big(\frac{d}{\nu}\Big)} {\Upsilon\Big(\frac{d}{\nu}\Big)}. \label{coisometry nu z-d eq 5}
\end{align}
Taking $z=0$ we get 
\begin{equation} \label{coisometry nu z-d eq 6}
\overline{\Upsilon(0)} = \frac{1} {\Upsilon\Big(\frac{d}{\nu}\Big)}.  
\end{equation}
Substituting this in equation (\ref{coisometry nu z-d eq 4}) we get
\begin{equation} \label{coisometry nu z-d eq 7*}
K_{(-d)}(-d) =\frac{1}{{\zeta_0}^2} \frac {\Upsilon\Big(\frac{d}{\nu}\Big)} {\Upsilon(0)}.
\end{equation}
From equations (\ref{coisometry nu z-d eq 3}), (\ref{coisometry nu z-d eq 5}) equating $\overline{\Upsilon(z)}$ and using equation (\ref{coisometry nu z-d eq 7*}), for all $z\in \mathbb{C}$  we get
\begin{align*} 
\frac{1}{{\zeta_0}^2 }  \frac{1}{{\zeta_0}^2 \Upsilon(0) K_{\Phi(z)}(-d)}=\frac{K_{z}\Big(\frac{d}{\nu}\Big)} {\Upsilon\Big(\frac{d}{\nu}\Big)}  \\ 
\frac{1}{{\zeta_0}^2 } \left(\frac{1}{{\zeta_0}^2 }   \frac{\Upsilon\Big(\frac{d}{\nu}\Big)}{\Upsilon(0)} \right)=K_{\Phi(z)}(-d) K_{z}\Big(\frac{d}{\nu}\Big) \\
\frac{1}{{\zeta_0}^2 } K_{(-d)}(-d)=K_{\Phi(z)}(-d) K_{z}\Big(\frac{d}{\nu}\Big). \\
\end{align*}
Taking $z=\Big(\frac{d}{\nu}\Big)$ we get 
\begin{align}
\frac{1}{{\zeta_0}^2 } K_{(-d)}(-d)&=\frac{1}{{\zeta_0}^2 } K_{(\frac{d}{\nu})}\Big(\frac{d}{\nu}\Big) \notag \\
\sum_{n=0}^{\infty} \frac{1}{{\zeta_n}^2} \, {{(\overline{-d})} }^n {(-d)}^{n} &=\sum_{n=0}^{\infty} \frac{1}{{\zeta_n}^2} \, {\left(\overline{\frac{d}{\nu} }\right) }^n {\Big(\frac{d}{\nu} \Big) }^{n} \notag  \\
\sum_{n=0}^{\infty} \frac{1}{{\zeta_n}^2} |d|^{2n} &=\sum_{n=0}^{\infty} \frac{1}{{\zeta_n}^2} \frac{|d|^{2n}}{|\nu|^{2n}}. \label{special equation for coisometry}
\end{align}
Case I. If $d=0$ then $\bar{\nu}d=0$ and by Lemma \ref{thm Entire fu space coisometry  Phi equal nu z} the required result follows because $K_{0}(z)=||K_{0}||^2=\frac{1}{{\zeta_0}^2}.$\\
Case II. If $d\neq 0,$ then comparing both sides of equation \eqref{special equation for coisometry} we get $|\nu|=1.$ Combining equations (\ref{coisometry nu z-d eq 4}), (\ref{coisometry nu z-d eq 5}) and (\ref{coisometry nu z-d eq 6}) we get
\begin{align}
\overline{\Upsilon(z)} &= {\zeta_0}^2 \overline{\Upsilon(0)}  K_{z}\Big(\frac{d}{\nu}\Big) \notag  \\
                     &= {\zeta_0}^2 \frac{1}{{\zeta_0}^2 \Upsilon(0) K_{(-d)}(-d)} K_{z}\Big(\frac{d}{\nu}\Big).  \label{coisometry nu z-d eq 6*}
\end{align}
Since $|\nu|=1,$ therefore, we have
 \begin{equation*} \label{coisometry nu z-d eq 7}
 K_{(-d)}(-d)=\sum_{n=0}^{\infty} \frac{1}{{\zeta_n}^2} \, {\bar{d}\,}^n {d\,}^{n}=\sum_{n=0}^{\infty} \frac{1}{{\zeta_n}^2} \, ({\nu}^n{\bar{d}\,}^n) ({\bar{\nu}}^n{d\,}^{n})=\sum_{n=0}^{\infty} \frac{1}{{\zeta_n}^2} \, {(\overline{\bar{\nu}d})}^n ({\bar{\nu}d})^{n}=K_{\bar{\nu}d}(\bar{\nu}d).
 \end{equation*}
Applying equation (\ref{norm of reproducing kernel}), for $p=\bar{\nu}d$ we get
\begin{equation}  \label{coisometry nu z-d eq 8}
K_{(-d)}(-d)=||K_{\bar{\nu}d}||^2.
\end{equation}
Substituting this in equation (\ref{coisometry nu z-d eq 6*}) we get $\overline{\Upsilon(z)}=\frac{1}{\Upsilon(0)} \frac {K_{z} (\bar{\nu} d)} { ||K_{\bar{\nu} d}||^2}$ because $\bar{\nu}=\frac{1}{\nu}.$ Taking conjugate we get 
\begin{align*}
\Upsilon(z)=\frac{1}{\overline{\Upsilon(0)}} \frac {K_{\bar{\nu} d} (z)} { ||K_{\bar{\nu} d}||^2} \,\, \,\,\, \text{for all}\, \, z \in \mathbb{C}.
\end{align*}
Combining equations (\ref{coisometry nu z-d eq 4}) and (\ref{coisometry nu z-d eq 8}) we get $|\Upsilon(0)|=\frac{1}{{\zeta_0} ||K_{\bar{\nu} d}||}.$ So, we can rewrite $$\Upsilon(z)=\alpha \frac {K_{\bar{\nu} d} (z)} { ||K_{\bar{\nu} d}||} \,\, \text{where} \,\,\alpha=\frac{1}{\overline{\Upsilon(0)} ||K_{\bar{\nu} d}||} \,\,  
\text{and}\,\, |\alpha|={\zeta_0}.$$\\
%%%%%%%%%%%%%%%%%%%%%%%%%%%%%%%%%%%%%%%%%%%%%%%%%%%%%%%%%%%%%%%%%%%%%%%%%

%\section*{Acknowledgements}

\textbf{Anuradha Gupta}\\
 Department of Mathematics,\\
  Delhi College of Arts and Commerce,\\
  University of Delhi, \\
  New Delhi-110023, India.\\
  \vspace{0.2cm}
 email: dishna2@yahoo.in\\
   \textbf{Geeta Yadav}\\
  Department of Mathematics,\\
   University of Delhi, \\
  New Delhi-110007, India.\\
  email: ageetayadav@gmail.com
\end{document}